\documentclass[11pt,reqno]{amsart}

\def\squarebox#1{\hbox to #1{\hfill\vbox to #1{\vfill}}}

\newcommand{\Ee}{{\mathcal E}}
\newcommand{\Mk}{\mu(\Ee)}

\newcommand{\w}{\omega}

\newcommand{\M}{{\mu(E)}}

\newcommand{\Z}{{\mathbb Z}}

\newcommand{\R}{{\mathbb R}}
\newcommand{\C}{{\mathbb C}}
\newcommand{\Ss}{{\bf S}}
\newcommand{\tm}{\widetilde{M}}

\newcommand{\N}{{\mathbb N}}
\newcommand{\T}{{\mathbb T}}
\newcommand{\I}{I_E}
\newcommand{\Ie}{I_{\E}}
\newcommand{\Me}{\mu(\Ee)}
\newcommand{\Tc}{\mu(\E)}
\newcommand{\Te}{\mathcal{T}_\E}
\newcommand{\U}{{U_E}}
\newcommand{\Ue}{{U_\E}}

\newcommand{\W}{{\mathcal{T_{\E}}}} 
\newcommand{\hM}{{\widehat{M}}}
\newcommand{\A}{{\mathcal {A}}}
\newcommand{\tM}{\widetilde{M}}
\newcommand{\ha}{{\widetilde{a}}}

\newcommand{\wi}{{{\widetilde{(\phi_k)_r}}}}
\newcommand{\wa}{{\widetilde{(f)_r}}}

\newcommand{\ov}{{\widehat{T}}}
\newcommand{\wideha}{{\widetilde{T}}}
\newcommand{\rr}{\mathbb R^+}
\newcommand{\Zz}{\mathbb Z^+}

\newcommand{\E}{{\bf E}}

\renewcommand{\Im}{\mathop{\rm Im}\nolimits}
\newcommand{\oM}{{\widehat{M}}}
\theoremstyle{plain}
\newtheorem{thm}{Theorem}

\newtheorem{lem}{Lemma}

\newtheorem{prop}{Proposition}

\newtheorem{deff}{Definition}

\theoremstyle{definition}

\title{Spectral results for operators commuting with translations on Banach spaces of sequences on $\Z^k$ and $\Z^+$}
\author[ V. Petkova]{ Violeta Petkova}
\address{Universit\'e Metz UFR MIM\\  Laboratoire de Mathématiques et Applications de Metz\\
                    UMR 7122\\
                    Ile du Saulcy 57045 Metz Cedex 1\\France}
\email{petkova@univ-metz.fr}


\numberwithin{equation}{section}
\begin{document}
\maketitle




\begin{abstract}
We study the spectrum of multipliers (bounded operators commuting with the shift operator $S$) 
on a Banach space $E$  of sequences on $\Z$. Given a multiplier $M$, 
we prove that $\overline{\tM(\sigma(S))}\subset \sigma(M)$ where $\tM$ is the symbol of $M$. 
We obtain a similar result for the spectrum of an operator commuting with the shift on a Banach space of sequences on $\Zz$.
 We generalize the results for multipliers on Banach spaces of sequences on $\Z^k$.
\end{abstract}

\section{Introduction}

Let $E\subset \C^\Z$ be a Banach space of complex sequences $(x(n))_{n \in \Z}.$ Denote by $S: \C^\Z \longrightarrow \C^\Z$,
 the shift operator defined by $Sx=(x({n-1}))_{n \in \Z}$, for $x=(x(n))_{n \in \Z} \in \C^\Z$, 
so that $S^{-1}x=(x({n+1}))_{n \in \Z}$. Let $F(\Z)$ be the set of sequences on $\Z$, which have a finite number of non-zero elements  and assume that $F(\Z)\subset E$. 
We will call a multiplier on $E$ every bounded operator $M$ on $E$ such that $MSa=SMa$, for every $a \in F(\Z)$. Denote by $\M$ the space of multipliers on $E$. For $z\in \T=\{z \in \C\::\:|z|=1\}$, consider the map $E \ni x \longrightarrow \psi_z(x)$ given by  $\psi_z(x)=(x(n)z^n)_{n \in \Z}$.
 Notice that if we assume that $\psi_z(E)\subset E$ 
for all $z\in \T$ and if for all $n \in \Z$, the map
 $$p_n: E \ni x\longrightarrow x(n)\in \C$$
 is continuous,
 then from the closed graph theorem it follows that the map $\psi_z$ is bounded on $E$. 
In this paper, we deal with Banach spaces of sequences on $\Z$ satisfying only the following three very natural hypothesis:\\
(H1) The set $F(\Z)$ is dense in $E$.\\
(H2) For every $n \in \Z$, $p_n$ is continuous from $E$ into $\C$. \\
(H3) We have $\psi_z(E)\subset E, \:\forall z \in \T$ and $\sup_{z \in \T}\|\psi_z\|<+\infty$.\\
We give some examples of spaces satisfying our hypothesis.\\
{\bf Example 1. } Let $\w$ be a weight (a sequence of positive real numbers) on $\Z$. Set 
$$l_\w^p(\Z)=\Bigl\{(x(n))_{n \in \Z}\in \C^\Z;\:\sum_{n \in \Z}|x(n)|^p\w(n)^p<+\infty\Bigr\}, \:1\leq p<+ \infty $$
and $\|x\|_{\w,p}=\Bigl(\sum_{n \in \Z}|x(n)|^p\w(n)^p\Bigr)^{\frac{1}{p}}.$
It is easy to see that the Banach space $l_\w^p(\Z)$ satisfies our hypothesis. Moreover, the operator $S$ (resp. $S^{-1}$)
 is bounded on $l_\w^p(\Z)$ if and only if 
$$\sup_{n\in \Z} \frac{\w(n+1)}{\w(n)}<+\infty\:\Big({\rm resp.} \:\sup_{n\in \Z} \frac{\w(n-1)}{\w(n)}<+\infty\Big).$$
{\bf Example 2.} Let ${\mathcal K}$ be a convex, non-decreasing, continuous function on $\R^+$ such that ${\mathcal K}(0)=0$ and
 ${\mathcal K}(x)>0$, for $x >0.$ For example, $\mathcal {K}$ may be $x^p$, for 
$\:1\leq p< +\infty$ or $x ^{p+\sin(\log(-\log( x)))}$, for $\: p>1+\sqrt{2}$. Let $\w$ be a weight on $\Z$. Set
$$l_{{\mathcal K},\w}(\Z)=\Bigl\{(x(n))_{n \in \Z}\in \C^\Z;\:\sum_{n \in \Z}{\mathcal K}\Bigl(\frac{|x(n)|}{t}\Bigr)\w(n)<+\infty, \:\rm {for\: some\: t>0}\Bigr\}$$
and
$\|x\|=\inf\Bigl\{ t>0,\:\sum_{n \in \Z}{\mathcal K}\Bigl(\frac{|x(n)|}{t}\Bigr)\w(n)\leq 1\Bigr\}.$
The space $l_{{\mathcal K},\w}(\Z)$, called a weighted Orlicz space (see \cite{R}, \cite{I}), is a Banach space satisfying our hypothesis. \\
{\bf Example 3.} Let $(q(n))_{n \in \Z}$ be a real sequence such that $q(n) \geq 1$, for all $n \in \Z$. 
For $a =(a(n))_{n \in \Z}\in \C^\Z$, set 
$\|a\|_{\{q\}}=\inf\Bigl\{ t>0,\:\sum_{n\in \Z}\Bigl|\frac{a(n)}{t}\Bigr|^{q(n)}\leq 1\Bigr\}.$
Consider the space $l^{\{q\}}=\{a\in \C^\Z;\:\|a\|_{\{q\}}<+\infty\},$ which is a Banach space (see \cite{Ed}) satisfying our hypothesis. Notice that if $\lim_{n \to +\infty}|q({n+1})-q(n)|\neq 0$ and if $\sup_{n \in \Z} q(n)<+\infty$, then either $S$ or $S^{-1}$ is not bounded (see \cite{N}). \\

It is easy to see that if $S(E)\subset E$, then by the closed graph theorem the restriction $S\vert_E$ of $S$ to $E$ is bounded from $E$ into $E$.
 From now on we will say that $S$ (resp. $S^{-1}$) is bounded when $S(E)\subset E$ (resp. $S^{-1}(E)\subset E$). 
If $S(E)\subset E$, we will call $\sigma(S)$ the spectrum of the operator $S$ with domain $E$. If $S$ is not bounded,
 denote by $\sigma(S)$ (resp. $\rho(S)$) the spectrum (resp. the spectral radius) of $\overline{S}$,
 where $\overline{S}$ is the smallest extension of $S\vert_{F(\Z)}$ as a closed operator. Recall that the domain  $D(\overline{S})$ of $\overline{S}$ is given by
$$D(\overline{S})=\{x\in E,\:\exists (x_n)_{n\in \N}\subset F(\Z)\:s.t. \:x_n\longrightarrow x\: and\: Sx_n\longrightarrow y\in E\}$$
and for $x \in D(\overline{S})$ we set $\overline{S}x=y$. 
We will denote by $\overset{\circ}{A}$ (resp. $\delta(A)$) the interior (resp. the boundary) of the set $A$. 
Denote by $e_k$ the sequence such that $e_k(n)=0$ (resp. 1), if $n\neq k$ (resp. $n=k$).
For a multiplier $M$ we set $\widehat{M}=M(e_0)$ and it is easy to see that
\begin{equation}\label{eq:M}
 Ma=\widehat{M}*a,\:\forall a\in F(\Z).
\end{equation}

Given $a\in \C^\Z$, define $\tilde{a}(z)=\sum_{n\in \Z} a(n)z^n$
and notice that if $a\in l^2(\Z)$, then $\tilde{a} \in L^2(\T).$
Denote by $\tm$ the function 
$$z\longrightarrow \widetilde{M}(z) = \sum_{n\in\Z}\hM(n)z^n.$$
Usually, $\widetilde{M}$ is called the {\bf symbol} of $M$.
It is easy to see that on the space of formal Laurent series we have the equality
\begin{equation}\label{eq:mult}
\widetilde{Ma}(z)=\tM(z)\ha(z),\:\:\forall z \in \C,\: \forall a \in F(\Z).
\end{equation}
However, it is difficult to determine for which $z\in \C$ the series $\tM(z)$ converges.
For $r>0$ let $C_r$ be the circle of center $0$ and radius $r$.
Recall the following result established in \cite{VZ} (Theorem 1).
\begin{thm} 
$1)$ If $S$ is not bounded, but $S^{-1}$ is bounded, then $\rho(S)=+\infty$ and if $S$ is bounded, but $S^{-1}$ is not bounded, then
$\rho(S^{-1})=+\infty.$\\
$2)$ We have $\sigma(S)=\Bigl\{ z\in \C,\:\frac{1}{\rho(S^{-1})}\leq |z|\leq \rho(S)\Bigr\}$.\\
$3)$ Let $M\in \M$. For $r> 0$ such that $C_r\subset \sigma(S)$, we have $\widetilde{M} \in L^\infty(C_r)$ and 
$|\widetilde{M}(z)|\leq \|M\|,$ a.e. on $C_r$.\\
$4)$ If $\rho(S)>\frac{1}{\rho(S^{-1})}$, then $\widetilde{M}$ is holomorphic on $\overset{\circ}{\sigma(S)}.$
\end{thm}
If $S$ and $S^{-1}$ are bounded, denote by $\I$ the interval $\Big [\frac{1}{\rho(S^{-1})}, \rho(S)\Big]$.
 If $S$ (resp. $S^{-1}$) is not bounded, denote by $\I$
 the interval $\Big [ \frac{1}{\rho(S^{-1})}, +\infty\Big[$ (resp. $]0, \rho(S)]$).  \\

The purpose of this paper is to use the symbol of an operator  $M \in \M$ in order to characterize its spectrum.
 We deal with three different setups. 
First we study the multipliers on $E$, next we examine Toeplitz operators on a Banach space of sequences on $\Zz$ 
and finally we deal with multipliers on a Banach space of sequences on $\Z^k.$ For $\phi\in F(\Z)$ we 
denote by $M_\phi$ the operator of convolution by $\phi$ on $E$. Let ${\mathfrak{S}}$ be
the closure with respect to the operator norm topology of the algebra generated by the operators $M_\phi$, 
for $\phi\in F(\Z).$ Our first result is
\begin{thm}
$1)$ If $M\in \M$, we have $\overline{\tM({\sigma(S)})}\subset \sigma(M)$.\\
$2)$ If $M\in {\mathfrak{S}}$ , then $\sigma(M)=\overline{\tM(\sigma(S))}.$
\end{thm}
Notice that here for a set $A$, we denote by $\overline{A}$ the closure of $A$. 
If $M\in \M$, then $\tM(\sigma(S))$ denotes the essential range of $\tM$ on $\sigma(S)$. 
Notice that $\tM$ is holomorphic on $\overset{\circ}{\sigma(S)}$ and essentially bounded on the boundary of $\sigma(S)$. 
In general we have no spectral calculus for the operators in $\M$ and it seems difficult to characterize the spectrum of $M \in \M$ 
without using its symbol.

We also study a similar spectral problem for Toeplitz operators. Let $\E\subset \C^{\Z^+}$ be a Banach space
and let $F(\Z^+)$ (resp. $F(\Z^-)$) be the space of the sequences on $\Z^+$ (resp. $\Z^-$)
 which have a finite number of non-zero elements. By convention, we will say that $x\in F(\Z)$ is a sequence of $F(\Z^+)$
 (resp. $F(\Z^-)$) if $x(n)=0$, for $n<0$ (resp. $n>0$). 
We will assume that $\E$ satisfies the following hypothesis:\\
(${\mathcal H}$1) The set ${F(\Z^+)}$ is dense in $ \E$.\\
(${\mathcal H}$2) For every $n \in \Z^+$, the application $p_n:x\longrightarrow x(n)$ is continuous from $\E$ into $\C$. \\
(${\mathcal H}$3) For $x=(x(n))_{n \in \Z+}\in \E$, we have $\gamma_z(x)=(z^n x(n))_{n \in \Z^+}\in \E$, 
for every $z \in \T$ and $\sup_{z\in\T}\|\gamma_z\|<+\infty.$\\
\begin{deff}
We define on $\C^{\Z^+}$ the operators $\Ss_1$ and $\Ss_{-1}$ as follows.
$$ For\: u\in \C^{\Z^+},\:(\Ss_1(u))(n)=0,\:if \: n = 0\:and \:(\Ss_1(u))(n)=u(n-1),\:if \:n\geq 1$$
$$(\Ss_{-1}(u))(n)=u(n+1),\: for\: n \geq 0.$$

\end{deff}
For simplicity, we note $\Ss$ instead of $\Ss_1$. It is easy to see that if $\Ss(\E)\subset \E$, then by the closed graph theorem the restriction $\Ss\vert_{\E}$ of $\Ss$ to $\E$ is bounded from $\E$ into $\E$. We will say that $\Ss$ (resp. $\Ss_{-1}$) is bounded when $\Ss(\E)\subset\E$ (resp $\Ss_{-1}(\E)\subset \E$). 
 Next, if $\Ss\vert_{\E}$ (resp. $\Ss_{-1}\vert_{\E})$ is bounded, $\sigma(\Ss)$ (resp. $\sigma(\Ss_{-1})$) denotes the spectrum of $\Ss\vert_{\E}$ (resp. $\Ss_{-1}\vert_{\E})$. 
 If $\Ss$ (resp. $\Ss_{-1}$) is not bounded, $\sigma(\Ss)$ (resp. $\sigma(\Ss_{-1})$) denotes the spectrum of the smallest closed extension of  $\Ss\vert_{F(\Z^+)}$ (resp. $\Ss_{-1}\vert_{F(\Z^+)}$).

For $u \in l^2(\Z^-)\oplus \E $ introduce
$$(P^+(u))(n)=u(n),\:\forall n\geq 0\:{\rm and}\: (P^+(u))(n)=0,\:\forall n<0.$$
If $S_{1}$ and $S_{-1}$ are the shift and the backward shift on $l^2(\Z^-)\oplus \E $, then $\Ss=P^+S_1$ and $\Ss_{-1}=P^+S_{-1}.$

{\bf Example 4.} Let $w$ be a positive sequence on $\Zz$.
Set 
$$l_w^p(\Zz)=\Bigl\{(x(n))_{n \in \Zz}\in \C^{\Zz};\:\sum_{n \in \Zz}|x(n)|^pw(n)^p<+\infty\Bigr\}, \:1\leq p<+ \infty $$
and $\|x\|_{w,p}=\Bigl(\sum_{n \in \Zz}|x(n)|^pw(n)^p\Bigr)^{\frac{1}{p}}.$
It is easy to see that the Banach space $l_w^p(\Zz)$ satisfies our hypothesis. The operator $\Ss$ (resp. $\Ss_{-1}$) is bounded on $l_w^p(\Zz)$, if and only if, $w$ satisfies
$$\sup_{n\in \Zz} \frac{w(n+1)}{w(n)}<+\infty \:\Big({\rm resp.} \:\sup_{n\in \Zz} \frac{w(n)}{w(n+1)}<+\infty\Big).$$

\begin{deff}
A bounded operator $T$ on $\E$ is called a Toeplitz operator, if we have:
$$(\Ss_{-1}T\Ss)u=Tu, \: \forall u \in F(\Z^+).$$
Denote by $\W$ the space of Toeplitz operators on $\E$.
\end{deff}
It is easy to see that if $T$ commutes either with $\Ss$ or with $\Ss_{-1}$, then $T$ is a Toeplitz operator.
Indeed, if $T\Ss_{-1} = \Ss_{-1} T$, then $T = \Ss_{-1} T \Ss.$
Notice that if $T\in \W$, we have
$Tu=P^+S_{-n}TS_n u,$ for all $u\in F(\Zz)$ and all $n>1$. Here $S_n$ (resp. $S_{-n}$) denotes $(S_1)^n$ (resp. $(S_{-1})^n$) where $S_1$ (resp. $S_{-1})
$ is the shift (resp. the backward shift) on $l^2(\Z^-)\oplus F(\Z^+)$.\\

Remark that we have $\Ss_{-1}\Ss=I$, 
however $\Ss\Ss_{-1} \not= I$ and this is the main difficulty in the analysis of Toeplitz operators.


Given a Toeplitz operator T, set $\widehat{T}(n)=(Te_0)(n)$ and $\widehat{T}({-n})=(Te_n)(0)$, 
for $n \geq 0$ and define $\widehat{T}=(\widehat{T}(n))_{n \in \Z}.$
It is easy to see that we have
\begin{equation}\label{eq:T} Tu=P^+(\widehat{T}*u),\:\forall u \in F(\Z^+).
\end{equation}


Set 
 $\wideha(z)=\sum_{n \in \Z}\widehat{T}(n)z^n,$
for $z \in \C$. Notice that the series $\widetilde{T}(z)$ could diverge.\\
If $\Ss$ and $\Ss_{-1}$ are bounded, we will denote by $\Ie$ the interval  $\Bigl[\frac{1}{\rho({\Ss}_{-1})}, \rho(\Ss)\Bigr]$. 
If $\Ss$ (resp. $\Ss_{-1}$) is not bounded, then $\Ie$ denotes
 $\Big[\frac{1}{\rho(\Ss_{-1})},+\infty\Big[$ \Big(resp. $\Big]0, {\rho(\Ss)}\Big]\Big)$. \\
If $\Ss$ and $\Ss_{-1}$ are bounded, denote by $\Ue$ the set $\Bigl\{ z \in \C,\:\frac{1}{\rho({\Ss}_{-1})}\leq |z|\leq \rho(\Ss)\Bigr\}$. 
If $\Ss$ (resp. $\Ss_{-1}$) is not bounded then $\Ue$ denotes $\Bigl\{ z \in \C,\:\frac{1}{\rho({\Ss}_{-1})}\leq |z|\Bigr\}$  $\Big($ resp. $ \Bigl\{z\in \C,\: |z|\leq \rho(\Ss)\Bigr\}\Big)$.
We have the following result (see Theorem 2 in \cite{VZ})
\begin{thm}
Let $T$ be a Toeplitz operator on $\E$.\\ 
$1)$ For $r \in \Bigl[\frac{1}{\rho({\Ss}_{-1})}, \rho(\Ss)\Bigr]$, if $\rho(\Ss)<+\infty$ or for $r \in \Bigl[\frac{1}{\rho({\Ss}_{-1})}, +\infty\Bigr[$, if $\rho(\Ss)=+\infty$ we have $\widetilde{T}\in L^\infty(C_r)$ and $|\widetilde{T}(z)|\leq \|T\|$, a.e. on $C_r$.\\
$2)$ If $\overset{\circ}{U_\E}$ is not empty, $\widetilde{T}\in \mathcal{H}^\infty(\overset{\circ}{U_\E}),$
where $\mathcal{H}^\infty(\overset{\circ}{U_\E})$ is the space of holomorphic and essentially bounded functions on $\overset{\circ}{U_\E}$.


\end{thm}
Denote by $\Tc$ the set of bounded operators on $\E$ commuting with either $\Ss$ or $\Ss_{-1}$.
 As mentioned above $\Tc\subset\Te.$ 
It is clear that the operators $(\Ss_n)_{n\geq 0}$ and $((\Ss_{-1})^n)_{n\geq 0}$ are included in $\Tc$.  
In this paper we prove the following

\begin{thm} If $\Ss$ and $\Ss_{-1}$ are bounded operators, we have 
\begin{equation} \label{eq:1.4}
\sigma (\Ss) = \{ z \in \C: |z| \leq \rho(\Ss)\}.
\end{equation}
\begin{equation}  \label{eq:1.5}
\sigma(\Ss_{-1}) = \{z \in \C:\: |z| \leq \rho(\Ss_{-1})\}.
\end{equation}
\end{thm}

 For the right $R$ and left $L$ weighted shifts on $l^2(\N)$ the results (\ref{eq:1.4}), (\ref{eq:1.5}) are classical (see for instance, \cite{R1}). 
Moreover, it is well known that the spectrum of $R$ and $L$ have a circular symmetry (\cite{R2}). The proofs of these results for $R$ and $L$ use
the structure of $l^2(\Z^+)$ and the analysis of the point spectrum is given by a direct calculus. 
In the general situation we deal with such an approach is not possible and our results on the symbols of Toeplitz operators play a crucial role. First we establish in Proposition 1 
the relation $\{\frac{1}{\rho(\Ss_{-1})} \leq |z| \leq \rho(\Ss)\} \subset \sigma(\Ss)$ and next we obtain (\ref{eq:1.4}). 
It seems that Theorem 4 is the first result concerning the description of $\sigma(\Ss)$ and 
$\sigma(\Ss_{-1})$ in the general setup when $(\mathcal H 1)$- $(\mathcal H 3)$ hold.\\

For operators commuting either with $\Ss$ or $\Ss_{-1}$ we have the following
\begin{thm} 
Suppose that $\Ss$ and $\Ss_{-1}$ are bounded.
Let $T$ be a bounded operator on $\E$ commuting with $\Ss$. 
Then we have 
\begin{equation}  \label{eq:1.6}
\overline{\widetilde{T}(\overset{\circ}{\sigma(\Ss)})}\subset \sigma(T).
\end{equation}
If $T$ is a bounded operator on $\E$ commuting with $\Ss_{-1}$, we have
\begin{equation} \label{eq:1.7}
\overline{\widetilde{T}(\overset{\circ}{\sigma(\Ss_{-1})})}\subset \sigma(T).
\end{equation}
\end{thm}
For $\phi\in \C^\Z,$ define
$$T_\phi f=P^+(\phi*f),\: \forall f\in \E.$$
If $\phi\in \C^{\Z}$ is such that $T_\phi$ is bounded on $\E$ (it is the case if for example $\phi \in F(\Zz)$), then $T_\phi\in \Tc$. 
The author has established similar results for multipliers and Winer-Hopf operators in weighted spaces $L_\w^2(\R)$ and $L_w^2(\rr)$
 (see \cite{V8}, \cite{Vc}). The spaces considered in this paper are much more general then weighted $l_\w^2(\Z)$ and $l_\w^2(\Zz)$ 
spaces.
 Here we consider not only Hilbert spaces, but also Banach spaces which may have a complicated structure (see Example 2 and Example 3).
 Moreover, we study multipliers on spaces where the shift is not a bounded operator. 
In these  general cases our spectral results are based heavily on the symbolic representation and this was the main motivation 
for proving the existence of symbols for the operators of the
classes we consider.

For $\phi\in C_c(\R^+)$, denote by $\mathcal{T}_\phi$ the operator defined on $L_\w^p(\R^+)$ by 
$$(\mathcal{T}_\phi f) (x)=P^+(\phi* f)(x),\:a.e .$$ 
Set $\beta_0=\lim_{t\to +\infty} \ln \|S_t\|^{\frac{1}{t}}$.
 A recent result of the author (see \cite{Vn}) shows that if $\mathcal{T}_\phi$ commutes with $(S_t)_{t\geq 0}$  
then 
$$\sigma(\mathcal{T}_\phi)=\hat{\phi}(V),$$
 where
$V=\{z\in \C,\: \Im z\leq\beta_0\}$.
 It is natural to conjecture that 
$$\sigma(T_\phi)=\hat{\phi}(U_\E)$$
 for $T_\phi$ with $\phi\in 
F(\Z)$ commuting with $\Ss$. \\

In Section 4, we study the so called joint spectrum for translation operators on a Banach space of sequences on $\Z^k$ 
and we generalize the results of Section 2.
In Theorem 7 we prove that the spectrum of a multiplier (bounded operator commuting with the translations) 
on a very general Banach space ${\mathcal E}$ of sequences on $\Z^k$ is related 
to the image under its symbol 
of the  joint spectrum of the translations $S_1$,...,$S_k$ (see Section 4 for the definitions).
 This joint spectrum denoted by $\Z_{\mathcal E}^k$ (see Section 4) is very important in our analysis. 
Notice that $\Z_{\mathcal E}^k \subset \sigma(S_1) \times ...\times \sigma(S_k)$ but
 in general the inclusion is strict. 
The fact that the symbol of a multiplier is holomorphic on the interior of $\Z_{\mathcal E}^k$ plays a crucial role.
  To our best knowledge it seems that Theorem 7 is the first result in the literature concerning the spectrum of operators commuting with translations on a Banach space of sequences on $\Z^k$.

\section{Spectrum of a multiplier}

First, consider the case of the multipliers on a Banach space $E$ satisfying (H1)-(H3) 
and suppose that $S$ or $S_{-1}$ is bounded on $E$.
Define $U_E=\{z\in \C,\: \frac{1}{\rho(S_{-1})}\leq |z|\leq \rho(S)\}$. 
To prove Theorem 2, we will need the following lemma established in \cite{VZ} (Lemma 4).
\begin{lem}
 For $\phi \in F(\Z)$, we have 
$|\widetilde{M_{\phi}}(z)|\leq \|M_\phi\|,\:\forall z\in \U.$

\end{lem}

\begin{deff}
For $a\in \C^\Z, $ and $r\in \R$, define the sequence $(a)_r$ so that
 $$(a)_r(n)=a(n)r^n,\:\forall n\in \Z.$$
\end{deff}

\begin{lem}
Let $r\in \I$ and $f\in E$ be such that $(f)_r\in l^2(\Z)$. 
If $M\in \M$, we have  $$(Mf)_r =(\oM_r*(f)_r),\:\:
 (\widetilde{Mf})_r(z)=\tM(rz)\widetilde{(f)_r}(z),\:\forall z\in \T$$ and $\widetilde{(Mf)_r}\in L^2(\T).$
\end{lem}
Lemma 2 is a generalization of (\ref{eq:M}).\\

{\bf Proof.}
The proof uses the arguments exposed in \cite{VZ} with some modifications. For the completeness we give here the details.
Let $M \in \M$. Let $(M_k)_{k \in \N}$ be a sequence such that 
$\lim_{k \to +\infty}\|M_k x-Mx\|=0 ,\: \forall x \in E,$
$\|M_k\|\leq \|M\|$  and $M_k=M_{\phi_k},$ where $\phi_k\in F(\Z),\:\forall k\in \N.$
The existence of this sequence is established in \cite{VZ} (Lemma 3).
Let $r\in I_E$. We have
$|\wi(z)|\leq \|M_{\phi_k}\|\leq \|M\|, \:\forall z \in \T,\:\forall k \in \N.$
We can extract from $\Bigl(\wi\Bigr)_{k \in \N}$ a subsequence which converges with respect to the weak topology $\sigma (L^{\infty}(\T), L^1(\T)) $ to a function $\nu_r \in L^\infty(\T).$ For simplicity, this subsequence will be denoted also by $\Big(\wi\Big)_{k \in \N}$. We obtain
$$\lim_{k \to +\infty}\int_\T \Big(\wi(z)g(z)-\nu_r(z)g(z)\Bigl)dz=0, \:\forall g \in L^1(\T)$$
and $\|\nu_r\|_{\infty}\leq \|M\|.$
Fix $f \in E$ such that $(f)_r\in l^2(\Z)$.
It is clear that 
$$\lim_{k \to +\infty}\int_\T \Big(\wi(z)\wa(z)g(z)-\nu_r(z)\wa(z)g(z)\Bigl)dz=0, \:\forall g \in L^2(\T).$$

We observe that the sequence $\Bigl( \wi\wa\Bigr)_{k \in \N}$ converges with respect to the weak topology of $L^2(\T)$ to $\nu_r \wa.$
Set $\widehat{\nu_r}(n)=\frac{1}{2\pi}\int_{-\pi}^{\pi}\nu_r(e^{it})e^{-itn}dt$, for $n \in \Z$ and let $\widehat{\nu_r}=(\widehat{\nu_r}(n))_{n \in \Z}$ be the sequence of the Fourier coefficients of $\nu_r$.
The Fourier transform from $l^2(\Z)$ to $L^2(\T)$ defined by
$\mathcal{F}:l^2(\Z)\ni (f(n))_{n \in \Z} \longrightarrow \tilde{f}\vert_\T \in L^2(\T)$
is unitary, so the sequence $\Bigl( (M_{\phi_k}f)_r \Bigr)_{k \in \N} =\Bigl((\phi_k)_r*(f)_r\Bigr)_{k \in \N}$ converges to $\widehat{\nu_r}*(f)_r$ with respect to the weak topology of $l^2(\Z)$. 
 Taking into account that $E$ satisfies (H2), for $n\in \Z$ we obtain 
$$\lim_{k \to +\infty}|((M_{\phi_k}f)_r-(Mf)_r)(n)|\leq\lim_{k \to +\infty}C\|M_{\phi_k}f-Mf\|=0.$$
Thus we deduce that 
$$(Mf)_r(n) =(\widehat{\nu_r}*(f)_r)(n),\: \forall n \in \Z,\: \forall f \in E ,$$
such that
$(f)_r\in l^2(\Z).$
This implies
$(\oM)_r*(f)_r=\widehat{\nu_r}*(f)_r, \forall f \in F(\Z)$
and then we get
$(\oM)_r=\widehat{\nu_r}.$
We conclude that 
$(Mf)_r =(\oM)_r*(f)_r,\:\forall f\in E$ such that $(f)_r\in l^2(\Z)$
and then we have 
$$(\widetilde{Mf})_r(z)=\tM(rz)\widetilde{(f)_r}(z),\:\forall z\in \T.$$
Since $\widetilde{(f)_r}\in L^2(\T)$ and $\tM\in L^\infty(U_E)$, it is clear that $\widetilde{(Mf)_r}\in L^2(\T).$
$\Box$\\

{\bf Proof of Theorem 2.}
Let $M\in \M$. Suppose that $\alpha\notin \sigma(M)$. Then we have $$(M-\alpha I)^{-1}\in \M.$$
 For $z\in \sigma(S)$, we obtain
$$\Big(\widetilde{(M-\alpha I)^{-1} f}\Big)(z)=
\Big(\sum_{n\in \Z} \widehat{(M-\alpha I)^{-1}}(n)z^n\Big)
 \Big(\sum_{n\in \Z} f(n)z^n\Big),$$
for all $ f\in E$, suc that for all $ r\in I_E, \:(f)_r\in l^2(\Z).$
If $g\in F(\Z)$, following Lemma 2, we may replace $f$ by $(M-\alpha I)g$.
 We get
$$\tilde{g}(z)=\Big(\sum_{n\in \Z} \widehat{(M-\alpha I)^{-1}}(n)z^n\Big)\Big( \sum_{n\in\Z} ((M-\alpha I)g)(n)z^n\Big)$$
$$=\widetilde{(M-\alpha I)^{-1}}(z)(\widetilde{Mg}(z)-\alpha \tilde{g}(z))
= \widetilde{(M-\alpha I)^{-1}}(z)(\tm(z)-\alpha)\tilde{g}(z),\:\forall g\in F(\Z),$$
for all $z\in \sigma(S)$.
This implies that for fixed $r\in \I$,
$$\widetilde{(M-\alpha I)^{-1}}(r\eta)(\tm(r\eta)-\alpha)=1,\:\forall \eta\in \T.$$
Since, $\widetilde{(M-\alpha I)^{-1}}$ is holomorphic on $\overset{\circ}{\sigma(S)}$ 
and essential bounded on $\delta(\sigma(S))$ (see Theorem 1), we obtain that 
$\tm(z)\neq \alpha$, for every $z\in \overset{\circ}{\sigma(S)}$ and for almost every $z\in \delta(\sigma(S))$. 
We conclude that 
$\tm({\sigma(S)})\subset \sigma(M),$
which proves the first part of the theorem. For te second one, let $M\in {\mathfrak{S}}.$ 
Then there exists a sequence $(M_{\phi_n})$ with $\phi_n\in F(\Z)$ such that 
$\lim_{n\to +\infty } \|M_{\phi_n}-M\|=0.$
Notice that from Lemma 1 it follows that
$$|\widetilde{M_{\phi_n}}(z)|\leq \|M_{\phi_n}\|\leq \|M\|, \:\forall z\in \U.$$
Taking into account that 
$|\widetilde{M_{\phi_n}}(z)-\widetilde{M_{\phi_k}}(z)|\leq \|M_{\phi_n}-M_{\phi_k}\|,\:\forall z\in \U,$
and the fact that $(M_{\phi_n})$ converges with respect to the norm operator theory, we conclude that 
$(\widetilde{M_{\phi_n}})$ converges uniformly on $\U$ to a function $\mu_M$. We observe that $(\widetilde{M_{\phi_n}})_r$ 
converges to
$\mu_M(r.)$ with respect to the weak topology $\sigma(L^\infty(\T), L^1(\T))$. So we can identify $\widetilde{M}(rz)$ and $\mu_M(rz)$ 
for $z\in \T$.

Consequently, $\tm$ is continuous on $\delta(\U)$. 
Let $\lambda \in \sigma(M)$. Then there exists a character $\gamma$ on ${\mathfrak{S}}$ such that $\lambda=\gamma(M)$. 
For $k\in \N^*$, denote by $S_k$ the operator $(S_1)^k$. 
We have 
$$\gamma(M_{\phi_n})=\gamma(\sum_{k\in \Z} \widehat{\phi_n}(k)S_k)=\sum_{k\in \Z} \widehat{\phi_n}(k) \gamma(S)^k$$

and we get
$\gamma(M)=\lim_{n\to +\infty} \gamma(M_{\phi_n})=\lim_{n\to+\infty} \widetilde{M_{\phi_n}}(\gamma(S))=\tm(\gamma(S)).$
We conclude that 
$\sigma(M)\subset \overline{\tm(\sigma(S))}.$
$\Box$

Now suppose that $M\notin {\mathfrak{S}}$. If $\alpha\in \sigma(M)$, then $\alpha=\gamma(M)$,
 where $\gamma$ is a character on the commutative Banach algebra $\M$.
 Following \cite{VZ} (Lemma 3), there exists a sequence $(M_{\phi_n})$, with $\phi_n\in F(\Z)$ such that
$\lim_{n \to +\infty}\|M_{\phi_n}a-Ma\|=0,\:\forall a\in E$
and we have $\lim_{n\to +\infty} \widetilde{M_{\phi_n}}(z)=\tM(z)$,
 for all $z\in \overset{\circ}{\sigma(S)}$ and for almost every $z\in \delta (\sigma(S))$.
If we suppose that $\gamma(S)\in \overset{\circ}{\sigma(S)}$ we have
$$\lim_{n\to +\infty} \gamma(M_{\phi_n})=\lim_{n\to +\infty} \widetilde{M_{\phi_n}}(\gamma(S))=\tM(\gamma(S)),$$
but in the general case we do not have
$$\lim_{n\to+\infty}\gamma(M_{\phi_n})=\gamma(M),$$
because $(M_{\phi_n})$ converges to $M$ with respect to the strong operator theory and may be not for the norm operator topology.

\section{Spectrum of an operator commuting either with $\Ss$ or $\Ss_{-1}$ on $\E$}
In this section we consider a Banach space $\E$ satisfying the conditions $({\mathcal{H}}_1)-({\mathcal{H}}_3)$.
 Suppose that $\Ss$ and $\Ss_{-1}$ are bounded on $\E$. 

Notice that it is easy to see that, for $\phi \in F(\Z)$, if $T_\phi$ commutes with $\Ss$ 
(resp. $\Ss_{-1}$) then $\phi\in F(\Z^+)$ (resp. $F(\Z^-)$). \\



\begin{lem}

  For $T\in \W$, $r\in \Ie$ and for $a\in \E$ such that $(a)_r\in l^2(\Z^+)$ we have\\
\begin{equation}{\label{eq:l0}}
(Ta)_r=P^+((\widehat{T})_r*(a)_r)
\end{equation}
and then $(Ta)_r\in l^2(\Z^+).$\\

\end{lem}
Lemma 3 is a generalization of the property (\ref{eq:T}).\\
{\bf Proof.} 
Let $T$ be a bounded operator in $\W$ and let $(\phi_k)_{k \in \N}\subset F(\Z)$ be such that
$$\lim_{k \to +\infty} \|T_{\phi_k}a-Ta\|=0, \:\forall a \in \E$$
and
$\|T_{\phi_k}\|\leq \|T\|, \:\forall k \in \N.$
The existence of the sequence $(T_{\phi_k})$ is established in \cite{VZ} (Lemma 5).
Fix $r \in \Ie$. 
We have (see Lemma 6 in \cite{VZ}),
$|\widetilde{(\phi_k)_r}(z)|\leq \|T_{\phi_k}\|\leq \|T\|, \:\forall z \in \T,\:\forall k \in \N.$
We can extract from $\Bigl(\widetilde{(\phi_k)_r}\Bigr)_{k \in \N}$ a subsequence which converges with respect to the weak topology
 $\sigma(L^{\infty}(\T), L^1(\T)) $ to a function $\nu_r \in L^\infty(\T).$ 
For simplicity, this subsequence will be denoted also by $\Big(\widetilde{(\phi_k)_r}\Big)_{k \in \N}$. 
Let $a\in \E$ be such that $(a)_r\in l^2(\Zz)$.
We conclude that,
 $\Bigl( \widetilde{(\phi_k)_r}\widetilde{(a)_r}\Bigr)_{k \in \N}$ converges with respect to the weak topology of
 $L^2(\T)$ to $\nu_r \widetilde{(a)_r}.$
Denote by $\widehat{\nu_r}=(\widehat{\nu_r}(n))_{n \in \Z}$ the sequence of the Fourier coefficients of $\nu_r$.
Since the Fourier transform from $l^2(\Z)$ to $L^2(\T)$ is an isometry, the sequence $(\phi_k)_r*(a)_r$ 
converges to $\widehat{\nu_r}*(a)_r$ with respect to the weak topology of $l^2(\Z)$. On the other hand,
 $\Big( T_{\phi_k}a \Big)_{k \in \N} $ converges to $Ta$ with respect to the topology of $\E$.
 Consequently, since $\E$ satisfies ($\mathcal{H}$2) we have 
$$\lim_{k \to +\infty}|((T_{\phi_k}a)_r-(Ta)_r)(n)| \leq\lim_{k \to +\infty}C\|T_{\phi_k}a-Ta\|=0,\:\forall n \in\N.$$
We conclude that 
\begin{equation}{\label{eq:l}}
(Ta)_r =P^+(\widehat{\nu_r}*(a)_r),\: \forall a \in \E\:{\rm such \: that}\:(a)_r\in l^2(\Zz).
\end{equation}
Since
$$(Ta)_r =P^+((\ov*a)_r),\: \forall a \in F(\Z^+),$$ it follows that
$\widehat{T}(n)r^n=\widehat{\nu_r}(n),\:\forall n \in \Z.$
Then (\ref{eq:l}) implies obviously (\ref{eq:l0}). 
Combining (\ref{eq:l0}) with the fact that $\widehat{T}\in l^\infty(\Z)$, 
it is clear that if $(a)_r\in l^2(\Zz)$, then $(Ta)_r\in l^2(\Zz)$. $\Box$\\


For the proof of Theorem 4 we need the following
\begin{prop}
Let $T$ be a bounded operator in $\mu(\E).$
Then we have 
$$\overline{\widetilde{T}(\overset{\circ}{\Ue})}\subset \sigma(T).$$
\end{prop}

{\bf Proof.} Let $T\in \Tc$ and suppose that $\lambda\notin \sigma(T)$. First we will show that $(T-\lambda I)^{-1}\in \Te.$ 
If $T\Ss = \Ss T$, then $(T- \lambda I) \Ss = \Ss (T - \lambda I)$ and we obtain $(T - \lambda I)^{-1}\Ss = \Ss (T - \lambda I)^{-1}$. 
As we have mentioned above this implies that $(T - \lambda I)^{-1}$ is a Toeplitz operator. 
In the same way we treat the case when $T\Ss_{-1} = \Ss_{-1}T.$\\
Set $h(n)=\widehat{(T-\lambda I)^{-1}} (n)$  
and fix $r\in \Ie$. 
For all $g\in \E$ such that $(g)_r\in l^2(\Zz)$, applying Lemma 3 with $(T-\lambda I)\in {\mathcal T}_{\E}$ and $a=g$ we get
  $(T - \lambda I)g\in l^2(\Zz)$. Then allpying a second time Lemma 3 with $(T-\lambda I)^{-1}\in {\mathcal T}_{\E}$ and 
$a=(T-\lambda I)g$, we get
$$(g)_r=P^+\Big((h)_r*((T-\lambda I)g)_r\Big),\:\forall g\in \E\:{\rm such\:\: that}\: (g)_r\in l^2(\Zz).$$
Since $(h)_r$ and $((T-\lambda I)g)_r$ are in $l^2(\Z^+)$ (see Lemma 3), we have 
$$
\|\widetilde{(g)_r}\|_{L^2(\T)}=\|(g)_r\|_{l^2(\Zz)}=\|P^+((h)_r*((T-\lambda I)g)_r)\|_{l^2(\Zz)}$$
$$\leq \|P^+\|\|(h)_r*((T-\lambda I)g)_r)\|_{l^2(\Zz)}=\|P^+\| \| \widetilde{(h)_r} 
(\widetilde{T}-\lambda) \widetilde{(g)_r}\|_{L^2(\T)}$$

$$\leq\|P^+\| \|\widetilde{(h)_r}\|_{L^\infty(\T)} 
 \|(\widetilde{T}-\lambda) \widetilde{(g)_r}\|_{L^2(\T)}$$

\begin{equation}{\label{eq:f}}
\leq C \|(\widetilde{T}-\lambda) \widetilde{(g)_r}\|_{L^2(\T)},
\:\forall g\in \E\:{\rm such\:that}\: \:(g)_r\in l^2(\Zz).
\end{equation}
First suppose that $1\in \overset{\circ}{\Ie}$. Then for $r=1$, we get
$$\|\widetilde{g}\|_{L^2(\T)}\leq C \|(\widetilde{T}-\lambda) \tilde{g}\|_{L^2(\T)},\:\forall g\in \E \cap l^2(\Zz).$$
 Assume that $\lambda=\widetilde{T}(z_0)$ for $z_0 \in \T \subset\overset{\circ}{U_{\E}}$.
According to Theorem 3, $\widetilde{T}$ is continuous on $\T$ and it is easy to
choose $f\in L^2(\T)$ so that
\begin{equation}\label{eq:in} 
2C\|(\widetilde{T}-\lambda)f\|_{L^2(\T)} < \|f\|_{L^2(\T)}.
\end{equation}
In fact, if $|\widetilde{T}(z) -\lambda| \leq \delta$ for $|z - z_0| < \eta(\delta) ,$ 
we take $f$ such that $f(z) = 0$ for $z$ s.t. $|z - z_0| \geq \eta(\delta)$ and $\|f\|_{L^2(\T)} = 1.$ 
For $\delta>0$ such that $2C \delta < 1$ we get the inequality (\ref{eq:in}).
Let $g\in l^2(\Z)$ be such that $f=\tilde{g}$ and
let $\beta=C\|\widetilde{T}-\lambda\|_\infty.$
Fix $\epsilon>0$ so that $\|\tilde{g}\|_{L^2(\T)}>(2\beta+2)\epsilon.$
Next let $g_\epsilon\in F(\Z)$ be such that 
$\|g_\epsilon-g\|_{l^2(\Z)}\leq \epsilon.$ Then we have
$$C\|(\widetilde{T}-\lambda)\widetilde{g_\epsilon}\|_{L^2(\T)}
\leq C\|(\widetilde{T}-\lambda)(\widetilde{g_\epsilon}-\tilde{g})\|_{L^2(\T)}
+C\|(\widetilde{T}-\lambda)\tilde{g}\|_{L^2(\T)}$$
$$\leq \beta \epsilon +\frac{1}{2}\|\tilde{g}\|_{L^2(\T)}
<\beta\epsilon+\frac{1}{2} \|\tilde{g}-\widetilde{g_\epsilon}\|_{L^2(\T)}+\frac{1}{2}\|\widetilde{g_\epsilon}\|_{L^2(\T)}
<(\beta+\frac{1}{2})\epsilon +\frac{1}{2}\|\widetilde{g_\epsilon}\|_{L^2(\T)}.$$
On the other hand, $\|\widetilde{g_\epsilon}\|_{L^2(\T)}\geq \|\tilde{g}\|_{L^2(\T)}-\epsilon \geq 2\beta +  \epsilon$, hence 
$(\beta +\frac{1}{2})\epsilon \leq \frac{1}{2}\|\widetilde{g_\epsilon}\|_{L^2(\T)}.$
This implies 
$$C\|(\widetilde{T}-\lambda)\widetilde{g_\epsilon}\|_{L^2(\T)}<\|\widetilde{g_\epsilon}\|_{L^2(\T)}.$$
Notice that for $f\in l^2(\Z)$ and $n\in \Zz$, we have 
$\widetilde{S_n f}(z)=z^n \tilde{f}(z), \:\forall z\in \T$.
Set $h=S_Ng_\epsilon$, where $N\in \Zz$ is chosen so  that $S_N g_\epsilon\in F(\Zz).$
We have 
$$C\|(\widetilde{T}-\lambda)\tilde{h}\|_{L^2(\T)}=C\|(\widetilde{T}-\lambda)\widetilde{S_N g_\epsilon}\|_{L^2(\T)}$$
$$=
C\|(\widetilde{T}-\lambda)\widetilde{g_\epsilon}\|_{L^2(\T)}<\|\widetilde{g_\epsilon}\|_{L^2(\T)}=\|\tilde{h}\|_{L^2(\T)}.$$
Taking into account (\ref{eq:f}), we obtain a contradiction and then $\widetilde{T}(z)\in \sigma(T)$ for $z\in \T$.\\
Now let $r\in\overset{\circ}{\Ie}$ and $r \neq 1$. Repeating the above argument, we choose  $g\in F(\Zz)$ so that 
$$C\|(\widetilde{T}-\lambda)\tilde{g}\|_{L^2(\T)}<\|\tilde{g}\|_{L^2(\T)}.$$
 Let $h$ be the sequence defined by $h(n)=g(n)r^{-n},\:\forall n\in \Zz$. Then $g=(h)_r$ and $h\in F(\Zz)$. 
We have 
$$C\|(\widetilde{T}-\lambda)\widetilde{(h)_r}\|_{L^2(\T)}<\|\widetilde{(h)_r}\|_{L^2(\T)}.$$
By using (\ref{eq:f}) once more, we obtain a contradiction and then $\widetilde{T}(C_r)\subset \sigma(T)$,
 where $C_r$ is the circle of center $0$ and radius $r$ and this completes the proof of the theorem. $\Box$\\

{\bf Proof of Theorem 4.} 
 The symbol of $\Ss$ is $z\longrightarrow z$ 
and according to Proposition 1, we have $U_{\E} \subset \sigma(\Ss).$ 
It remains to show that $\{ z\in \C,\:|z| < \frac{1}{\rho(\Ss_{-1})}\} \subset \sigma(\Ss).$
 We apply the argument of \cite{Vc}.\\
 For $0 < |z| < \frac{1}{\rho(\Ss_{-1})}$ we write
\begin{equation} \label{eq:ss}
\Ss_{-1} - \frac{1}{z} = - \frac{1}{z}\Bigl(\Ss_{-1} (\Ss - z)\Bigr).
\end{equation}
If $z \notin \sigma(\Ss)$, then there exists $g \not= 0$ such that $(\Ss - z)g = e_0.$ This implies $(\Ss_{-1} - \frac{1}{z}) g = 0$ and we obtain a contradiction with the fact that $\frac{1}{|z|} > \rho(\Ss_{-1}).$ 
This completes the proof of (\ref{eq:1.4}). 

Now we pass to the analysis of $\sigma(\Ss_{-1}).$ As above assume that $\Ss_{-1}$ is bounded. Following \cite{Vc}, we show first 
that for the approximative spectrum $\Pi(\Ss)$ of $\Ss$ we have
$$\Pi(\Ss) \subset \{z\in \C, \:\frac{1}{\rho(\Ss_{-1})} \leq |z| \leq \rho(\Ss)\}.$$
In fact, for $z \not= 0,$ if there exists a sequence $f_n,\: \|f_n\| = 1$ such that $(\Ss - z) f_n \to 0$ as $n \to \infty,$ then from
(\ref{eq:ss}) 
we deduce that $(\Ss_{-1} - \frac{1}{z})f_n \to 0$ and this yields $\frac{1}{|z|} \leq \rho(\Ss_{-1}).$
On the other hand, if $0 \in \Pi(\Ss)$, there exists a sequence $f_n,\: \|f_n\| = 1$ such that $\Ss f_n \to 0$ and this yields a contradiction with the equality $f_n = \Ss_{-1}\Ss f_n.$\\

For the proof of (\ref{eq:1.5}) we use for $z \not= 0$ the adjoint operators $\Ss^*, \: \Ss_{-1}^*$ and the equality
$$z\Bigl(\frac{1}{z} I -\Ss^*) = \Ss^* \Bigl((\Ss_{-1})^* - z I 
\Bigr).$$
The symbol of $\Ss_{-1}$ is $z\longrightarrow \frac{1}{z}$ and an application of Proposition 1 yields 
$\{\frac{1}{\rho(\Ss)} \leq |z| \leq \rho(\Ss_{-1})\} \subset \sigma(\Ss_{-1})$.
 Next assume that $0 < |z| < \frac{1}{\rho(\Ss)}$. 
We are going to repeat  the argument of the proof of Theorem 3 in \cite{Vc} and  for completeness we give the proof.
 First, $0 \in \sigma_r(\Ss)$, $\sigma_r(\Ss)$ being the residual spectrum of $\Ss$.
 In fact, if this is not true, 0 will be in $\Pi(\Ss)$ and this is a contradiction. 
Secondly, we deduce that 0 will be an eigenvalue of the adjoint operator $\Ss^*$.
Let $\Ss^* g = 0$ with $g \not= 0.$ 
If $(\Ss_{-1})^*- z I $ is surjective, than there exists $f \not= 0$ such that $( (\Ss_{-1})^* - z)f = g$ and we get
 $(\frac{1}{z} - \Ss^*)f = 0$, hence $\frac{1}{|z|} \leq \rho(\Ss^*)= \rho(\Ss)$ which is impossible. 
Thus $z \in \sigma(\Ss_{-1}^*)$ and, passing to the adjoint, we complete the proof. $\Box.$

For the proof of Theorem 5 we need the following

\begin{lem}
Let $\phi\in F(\Z^+)$ $($resp. $F(\Z^-))$. Then for $z \in \sigma(\Ss)$ (resp. $z\in \sigma(\Ss_{-1})$), we have 
$$ |\widetilde{(\phi)}(z)|\leq \|T_{\phi}\|\leq \|T\|. $$
\end{lem}

{\bf Proof.} 
Suppose that $|z| = \rho(\Ss).$ Then $z$ is in $\Pi(\Ss)$ and there exists a sequence
 $(f_n)_{n \in \N} \subset \E$ such that $\|f_n\|=1$ and $\lim_{n\to +\infty} \|\Ss f_n-zf_n\|=0$. 
Then  for $\phi\in F(\Z^+)$, we have for some $N>0$,
$$\|\phi*f_n-\widetilde{\phi}(z)f_n\| \leq \sum_{k=0}^{N}(\sup_{|k|\leq N}|\phi(k)|) \|\Ss^k f_n-z^kf_n\|$$ and we obtain
$$\lim_{n \to +\infty}\|\phi*f_n-\widetilde{\phi}(z)f_n\|=0.$$
Since $$|\widetilde{\phi}(z)|=\|\widetilde{\phi}(z)f_n\|\leq \|\widetilde{\phi}(z)f_n-\phi*f_n\|+\|T_\phi f_n\|,$$ it
 follows that $|\widetilde{\phi}(z)|\leq \|T_\phi\|.$ By using the maximum principle for the analytic function $\widetilde{\phi}$ 
we complete the proof for $z\ \in \sigma(\Ss).$ For $\sigma(\Ss_{-1})$ we apply the same argument. $\Box.$

\begin{lem} Let $T$ be a bounded operator on $\E$ commuting with $\Ss$ (resp. $\Ss_{-1}$). Let
 $r$ be such that there exists $z\in \sigma(\Ss)$ (resp. $\sigma(\Ss_{-1})$) with $r=|z|$. Then  for $a\in \E$ such that
 $(a)_r\in l^2(\Z^+)$ (resp. $(a)_r \in l^2(\Z^-)$) we have
\begin{equation}{\label{eq:l0}}
(Ta)_r=P^+((\widehat{T})_r*(a)_r)
\end{equation}
and then $(Ta)_r\in l^2(\Z^+)$ (resp. $(Ta)_r \in l^2(\Z^-)$). 

\end{lem}

{\bf Proof.} For the  proof we apply Lemma 4 and the same arguments as those in the proof of Lemma 3. $\Box$\\

By using Lemmas 4-5 and repeating the arguments of the proof of Proposition 1, we obtain Theorem 5. 

We leave the details to the reader.\\

\section{Spectral results for multipliers on Banach space of functions on $\Z^k$}

Let $F(\Z^k)$ be the space of sequences of $\Z^k$ with a finite number of not vanishing  terms. Let $\Ee$ be a Banach space of sequences on $\Z^k$ satisfying the following conditions:\\
$(h_1)$ $F(\Z^k)$ is dense in $\Ee$.\\
$(h_2)$ For every $n\in \Z^k$, the application $\Ee\ni x\longrightarrow x(n)\in \C^k$ is continuous.\\
$(h_3)$ For every $z\in \T^k$, we have $\psi_z(\Ee)\subset \Ee$ and $\sup_{z\in \T^k}\|\psi_z\|<+\infty,$ where
$$(\psi_z(x))(n_1,...,n_k)=x(n_1,...,n_k)z_1^{n_1}...z_k^{n_k},\:\forall n\in \Z^k, \:x\in \Ee.$$


Denote by $\Mk$ the space of bounded operators on $\Ee$ commuting with the translations.
Denote by ${\mathcal{S}}_i$ the operator of translation by $e_i$, where 
$e_i(n)=1,$ if $n_i=1$ and $n_j=0$, for $j\neq i$ and else $e_i(n)=0$.
 Suppose that the operator ${\mathcal{S}}_i$ is bounded on $\Ee$ for all $i\in \Z$.
For $M\in \mu(\Ee)$ define
$$\widetilde{M}(z)=\sum_{n\in \Z^k} \widehat{M}(n_1,...,n_k)z_1^{n_1}...z_k^{n_k},$$
for $z=(z_1,...,z_k)\in \C^k$, where
$\widehat{M}(n_1,...,n_k)=M(e_0)(n_1,...,n_k).$
For $a\in \Ee$, set 
$$\tilde{a}(z)=\sum_{n\in \Z^k} a(n)z_1^{n_1}...z_k^{n_k},\:\forall z\in \C^k.$$
Notice that for $M\in \Me$ and $a\in F(\Z^k)$, we have
$$Ma=\widehat{M}* a, \:\forall a\in F(\Z^k)$$

and formally we get 
$$\widetilde{Ma}(z)=\widetilde{M}(z)\tilde{a}(z), \: a\in \Ee,\: z\in \C^k.$$

If $\phi\in F(\Z^k)$ denote by $M_\phi$ the operator given by
$M_\phi f=\phi*f,\:\forall f\in \Ee.$ Define the set 
$$\Z^k_{\Ee}=\Big\{z\in \C^k,\:\Big|\sum_{n\in \Z^k} \phi(n)z_1^{n_1}..z_k^{n_k}\Big|\leq \|M_\phi\|, \:\forall \phi \in F(\Z^k)\Big\}.$$

Denote by $\sigma_A(B_1,..., B_p)$ the joint spectrum of the elements $B_1$,...,$B_p$ in a commutative Banach algebra $A$. 
Recall that $\sigma_A(B_1,...,B_p)$ is the set of $(\lambda_1,...,\lambda_p) \in \C^p$ such that for all $L\in A$, the operator 
$(B_1-\lambda_1 I)L+...+(B_p-\lambda_p I) L$ is not invertible  (see \cite{Z}). 
We have also the representation
$$\sigma_A(B_1,..., B_p)=\{(\gamma(B_1),...,\gamma(B_p)): \: \gamma  \:{\text is \:a\: character\: on}\: A\}.$$
 It is clear that 
$$ \sigma_A(B_1,...,B_p)\subset\sigma(B_1)\times...\times \sigma(B_p), $$
 but in general these two sets are not equal and the inclusion could be strict.
\begin{deff}
Denote by ${\mathcal{A}}$ the closure of the subalgebra generated by the operators ${M_\phi}$, $\phi\in F(\Z^k)$, with respect
 to the operator norm topology.
\end{deff}

\begin{prop}
 We have $\sigma_{\A}({\mathcal{S}}_1,...., {\mathcal{S}}_k)=\Z_{\Ee}^k$.
\end{prop}
{\bf Proof.} Let $z\in \C^k$ be such that 
$z=(\gamma({\mathcal{S}}_1),..., \gamma({\mathcal{S}}_k)),$
where $\gamma$ is a character on the algebra ${\mathcal{A}}$. Then 
$$\sum_{n\in \Z^k} \phi(n)z_1^{n_1}...z_k^{n_k}=\sum_{n\in \Z^k}
 \phi(n)\gamma({\mathcal{S}}_1)^{n_1}...\gamma({\mathcal{S}}_k)^{n_k}=\gamma(M_\phi),\forall\phi \in F(\Z^k)$$
and it is clear that 
$|\gamma(M_\phi)|\leq \| M_\phi\|,\:\forall \phi \in F(\Z^k).$
It follows that $\sigma_{\mathcal{A}}({\mathcal{S}}_1,...., {\mathcal{S}}_k)\subset \Z_{\Ee}^k.$
On the other hand, if $z\in \Z_{\Ee}^k$, we define
$$\gamma_z: M_\phi\longrightarrow \sum_{n\in \Z^k} \phi(n)z_1^{n_1}...z_k^{n_k}.$$
The application $\gamma_z$ is a character on ${\mathcal{A}} $
 and this implies that $z=(\gamma_z({\mathcal{S}}_1),...,\gamma_z({\mathcal{S}}_k))$ is in the joint spectrum of 
${\mathcal{S}}_1,..., {\mathcal{S}}_k$ in ${\mathcal {A}}$. 
 So we have $ \Z_{\Ee}^k=\sigma_{\mathcal{A}}({\mathcal{S}}_1,...., {\mathcal{S}}_k).$ $\Box$

Define
$$I_{\Ee}=\{r\in \R^k, \:r_1\T\times...\times r_k\T \in \overset{\circ}{\Z^k_{\Ee}} \:\}.$$ 
For $a\in \Ee$ and $r\in\C^k$, denote by $(a)_r$ the sequence
$$(a)_r(n_1,...,n_k)=a(n_1,...,n_k)r_1^{n_1}...r_k^{n_k},\:\forall (n_1,...,n_k)\in \Z^k.$$

The following theorem was established in \cite{VL} (Theorem 4 and Collorary 1).
\begin{thm}
 Let $\Ee$ be a Banach space of sequences on $\Z^k$ satisfying ($h_1$), ($h_2$) and ($h_3$) 
and such that ${\mathcal{S}}_i$ is bounded on $\Ee$ for all $i\in \Z$. Suppose that $\overset{\circ}{\Z^k_{\Ee}}\neq\emptyset $. 
Then, for $M\in \Mk$, there exists $ \theta_M \in {\mathcal{H}}^\infty(\overset{\circ}{\Z^k_\Ee})$ such that for $f\in {F}(\Z^k)$
 we have
$\widetilde{Mf}(z)=\theta_M(z)\tilde{f}(z),\:\: \forall  z \in \overset{\circ}{\Z^k_{\Ee}}.$
\end{thm}

Following (Lemma 2 in \cite{VL}) there exists a sequence $(M_m)_{m\in \N}\subset \Me$ such that:\\
$\lim_{m\to +\infty}\| M_m a - Ma\|=0,\:\forall a\in \Ee,$
$M_m=M_{\phi_m}$, where $\phi_m\in F(\Z^k)$ and $\|M_m\|\leq C \|M\|$. 
Notice that using the sequence $(M_m)$ and the same arguments as in the proof of Lemma 2, we obtain that in fact 
$\theta_M=\widetilde{M}$  and, moreover, we get the following

\begin{lem}
 For $M\in \Me$ and for $f\in \Ee$ such that $(f)_r\in l^2(\Z^k)$, for all $r\in I_{\Ee}$ we have
$$\widetilde{Mf}(z)=\widetilde{M}(z)\tilde{f}(z),\:\: \forall  z \in \overset{\circ}{\Z^k_{\Ee}}.$$ 

\end{lem}


Now we obtain the following spectral result.
\begin{thm}
1) For $M\in \Me$, we have $\overline{\widetilde{M}(\overset{\circ}{\Z^k_{\Ee})}}\subset \sigma(M).$\\
2) For $M\in \mathcal{A}$, we have $\tM(\Z_\Ee^k)=\sigma(M).$
 
\end{thm}


{\bf Proof.} Let $M\in \Me$.
Suppose that $\alpha\notin \sigma(M)$. Then we have $K=(M-\alpha I)^{-1}\in \Me$ and 
\begin{equation}\label{eq:k}
\widetilde{Kf}(z)=\widetilde{K}(z)\tilde{f}(z),\:\forall z\in \overset{\circ}{\Z^k_{\Ee}},
\:\forall f\in \Ee\:s. t. \: (f)_r\in l^2(\Z^k), 
\:\forall r\in I_E.
\end{equation}
$$\Big(\widetilde{(M-\alpha I)^{-1} f}\Big)(z)
=\Big(\sum_{n\in \Z} \widehat{(M-\alpha I)^{-1}}(n)z^n\Big) \Big(\sum_{n\in \Z} f(n)z^n\Big),\:\forall f\in \Ee,
 s.t. \:\forall r\in I_{\Ee}, \:(f)_r\in l^2(\Z^k).$$
If $g\in F(\Z^k)$, following Lemma 6, we may replace  $f$ by $(M-\alpha I)g$ in (\ref{eq:k}).
 We get
$$\tilde{g}(z)
=\Big(\sum_{n\in \Z} \widehat{(M-\alpha I)^{-1}}(n)z^n\Big)\Big( \sum_{n\in\Z} ((M-\alpha I)g)(n)z^n\Big)$$
$$=\widetilde{K}(z)(\widetilde{Mg}(z)-\alpha \tilde{g}(z))= \widetilde{K}(z)(\tm(z)-\alpha)\tilde{g}(z),$$
for all $z  \in \overset{\circ}{\Z^k_{\Ee}}$.
This implies that for fixed $r\in I_{\Ee}$,
$\widetilde{K}(r\eta)(\tm(r\eta)-\alpha)=1,\:\forall \eta\in \T^k.$
Since, $\widetilde{K}$ is holomorpic on $\overset{\circ}{\Z^k_{\Ee}}$, we obtain that 
$\tm(z)\neq \alpha$, for every $z\in \overset{\circ}{\Z^k_{\Ee}}$. 
We conclude that 
$$\tm(\overset{\circ}{\Z^k_{\Ee}})\subset \sigma(M),$$
which proves part 1).
Now suppose that $M=M_\phi$, with $\phi\in F(\Z^k).$ Let $\lambda\in \sigma(M_\phi).$
Then there exists $\gamma$ a character on $\mu(\Ee)$ such that 
$$\lambda=\gamma(M_\phi)=\sum_{n\in \Z^k}\phi(n)\gamma ({\mathcal{S}}_{n_1,...,n_k})
=\tilde{\phi}(\gamma({\mathcal{S}}_1),...,\gamma({\mathcal{S}}_k))\in \tilde{\phi}(\Z_\Ee^k).$$
The end of the proof of 2) is now very similar to the proof of 2) in Theorem 2 and is left to the reader.   
$\Box$

\end{document}